\documentclass[11pt]{amsart}

\usepackage{amsthm,amssymb,amsfonts,amsmath,enumerate,graphicx}

\newcommand{\Z}{\mathbb{Z}}
\newcommand{\R}{\mathbb{R}}
\newcommand{\K}{{\mathcal K}}

\renewcommand{\P}{{\mathcal P}}
\newcommand{\T}{{\mathcal T}}
\def\th{^{\text{th}}}
\newcommand\fl[1]{\left\lfloor {#1} \right\rfloor} 
\newcommand\fr[1]{\left\{ {#1} \right\}} 
\newcommand\saw[1]{\left(\!\left( #1 \right)\!\right)}
\newcommand\ded{{\rm s}}
\newcommand\car{{\rm c}}
\newcommand\drc{{\rm drc}}
\def\a{{\boldsymbol a}}
\def\e{{\boldsymbol e}}
\def\m{{\boldsymbol m}}
\def\u{{\boldsymbol u}}
\def\v{{\boldsymbol v}}
\def\z{{\boldsymbol z}}

\newtheorem{theorem}{Theorem}
\newtheorem{corollary}[theorem]{Corollary}

\newtheorem{lemma}[theorem]{Lemma}

\newtheorem*{definition}{Definition}

\newcommand\comment[1]{}                
\renewcommand\comment[1]{\emph{[#1]}}           

\title{Dedekind--Carlitz Polynomials as Lattice-Point Enumerators in Rational Polyhedra}

\author{Matthias Beck}
\address{Department of Mathematics\\
         San Francisco State University\\
         San Francisco, CA 94132, USA}
\email{beck@math.sfsu.edu}
\urladdr{http://math.sfsu.edu/beck}

\author{Christian Haase}
\address{Fachbereich Mathematik \& Informatik \\
  Freie Universit\"at Berlin \\
  14195 Berlin\\
  Germany}
\email{christian.haase@math.fu-berlin.de}
\urladdr{http://ehrhart.math.fu-berlin.de}
\thanks{Research of Haase supported by DFG Emmy Noether fellowship HA 4383/1. 
We thank Robin Chapman, Eric Mortenson, and an anonymous referee for helpful comments.} 

\author{Asia R. Matthews}
\address{Department of Mathematics \& Statistics\\
Queen's University\\
Kingston, ON, K7L 3N6, Canada}
\email{asiamath@mast.queensu.ca}



\date{8 February 2008}

\begin{document}

\begin{abstract}
We study higher-dimensional analogs of the \emph{Dedekind--Carlitz polynomials}
\[
  \car \left( u, v; a, b \right) := \sum_{ k=1 }^{ b-1 } u^{ \fl{ \frac{ ka }{ b } } } v^{ k-1 } ,
\]
where $u$ and $v$ are indeterminates and $a$ and $b$ are positive integers. Carlitz proved that these polynomials satisfy the \emph{reciprocity law}
\[
  \left( v-1 \right) \, \car \left( u,v; a,b \right) +  \left( u-1 \right) \, \car \left( v,u;b,a \right) = u^{ a-1 } v^{ b-1 } -1\, ,
\]
from which one easily deduces many classical reciprocity theorems for the Dedekind sum and its generalizations.
We illustrate that Dedekind--Carlitz polynomials appear naturally in generating functions of rational cones and use this fact to give geometric proofs of the Carlitz reciprocity law and various extensions of it. Our approach gives rise to new reciprocity theorems and computational complexity results for Dedekind--Carlitz polynomials, a characterization of Dedekind--Carlitz polynomials in terms of generating functions of lattice points in triangles, and a multivariate generalization of the Mordell--Pommersheim theorem on the appearance of Dedekind sums in Ehrhart polynomials of 3-dimensional lattice polytopes.
\end{abstract}

\maketitle


\section{Introduction}\label{introsection}

While studying the transformation properties of   
$\eta (z) := e^{ \pi i z / 12 } \prod_{ n \geq 1 } \left( 1 - e^{ 2 \pi i n z } \right) $
under $ \mbox{SL}_{2} ( \Z ) $, Richard Dedekind, in the 1880's \cite{dedekind}, naturally arrived at what we today call the \emph{Dedekind sum}
\[
  \ded \left( a, b \right) := \sum_{ k=0 }^{ b-1 } \saw{ \frac{ ka }{ b } } \saw{ \frac{ k }{ b } } ,
\]
where $a$ and $b$ are positive integers, and $\saw x$ is the \emph{sawtooth function} defined by
\[
  \saw x :=
  \begin{cases}
    \fr x - \frac 1 2 & \text{ if } x \notin \Z , \\
    0 & \text{ if } x \in \Z .
  \end{cases}
\]
Here $\fr x = x - \fl x$ denotes the \emph{fractional part} of $x$.
The Dedekind sum and its generalizations have since intrigued mathematicians from various areas
such as analytic \cite{almkvist,dieterdedekind} and algebraic number theory \cite{meyerdedekind,solomondedsum}, topology \cite{hirzebruchzagier,meyersczech}, 
algebraic \cite{brion,pommersheim} and combinatorial geometry \cite{ccd,mordell}, and algorithmic complexity \cite{knuth}. 

Almost a century after the appearance of Dedekind sums, Leonard Carlitz introduced the following polynomial generalization, which we will call a \emph{Dedekind--Carlitz polynomial}:
\[
  \car \left( u, v; a, b \right) := \sum_{ k=1 }^{ b-1 } u^{ \fl{ \frac{ ka }{ b } } } v^{ k-1 } .
\]
Here $u$ and $v$ are indeterminates and $a$ and $b$ are positive integers.
Undoubtedly the most important basic property for any Dedekind-like sum is \emph{reciprocity}. For the Dedekind--Carlitz polynomials, this takes on the following form \cite{carlitzpolynomials}.

\begin{theorem}[Carlitz]\label{carlitzreciprocity}
If $u$ and $v$ are indeterminates, and $a$ and $b$ are relatively prime positive integers, then
\[
  \left( v-1 \right) \, \car \left( u,v; a,b \right) +  \left( u-1 \right) \, \car \left( v,u;b,a \right) = u^{ a-1 } v^{ b-1 } -1\, .
\]
\end{theorem}

Carlitz's reciprocity theorem generalizes that of Dedekind \cite{dedekind}, which states that for relatively prime positive integers $a$ and $b$,
\begin{equation}\label{dedreci}
  \ded \left( a,b \right) + \ded \left( b,a \right) = - \frac{ 1 }{ 4 }+ \frac{ 1 }{ 12 } \left( \frac{ a }{ b } + \frac{ 1 }{ ab } + \frac{ b }{ a } \right) .
\end{equation}
Dedekind reciprocity follows from Theorem~\ref{carlitzreciprocity} by applying the operators $u \, \partial{u}$ twice and $v\, \partial v$ once to Carlitz's reciprocity identity and converting the greatest-integer functions into fractional parts.

Our motivation for this paper stems from the appearance of Dedekind sums in lattice-point enumerators for rational polyhedra.
The first such instance was discovered by Mordell \cite{mordell} (the case $t=1$ in the following theorem) and vastly generalized by Pommersheim \cite{pommersheim}, who laid the foundation for the appearance of Dedekind-like sums in \emph{Ehrhart polynomials} $L_\P (t) := \# \left( t \P \cap \Z^d \right)$; here $\P$ is a lattice $d$-polytope (i.e., $\P$ has integral vertices) and $t$ denotes a positive integer.

\begin{theorem}[Mordell--Pommersheim]\label{thm.L}
Let $\T$ be the convex hull of $(a,0,0)$, $(0,b,0)$, $(0,0,c)$, and $(0,0,0)$, where $a, b$ and $c$ are pairwise relatively prime positive integers. Then the Ehrhart polynomial of $\T$ is
\begin{align*}
  L_\T (t) = &\frac{abc}{6} \, t^3 + \frac{ab+ac+bc+1}{4} \, t^2 \\
  &+ \left( \frac{3}{4} + \frac{a+b+c}{4} + \frac{1}{12} \left( \frac{bc}{a} + \frac{ca}{b} + \frac{ab}{c} + \frac{1}{abc} \right) - \ded \left( bc, a \right) - \ded \left( ca, b \right) - \ded \left( ab, c \right) \right) t + 1 \, .
\end{align*}
\end{theorem}

There are natural geometric interpretations for three of the four coefficients appearing in the Ehrhart polynomial of the Mordell--Pommersheim tetrahedron $\T$: 
the leading coefficient is the volume of $\T$, the second leading coefficient equals half of the sum of the areas of the faces of $\T$, and the constant term is the Euler characteristic of $\T$. 
(In fact, similar interpretations hold for any $d$-dimensional lattice polytope \cite{ccd,ehrhartpolynomial}.)
However, aside from toric varieties attached to $\T$ \cite{pommersheim}, a geometric reason for the appearance of Dedekind sums in the linear term of the Ehrhart polynomial in Theorem~\ref{thm.L} has so far eluded mathematicians.  In this paper we attempt to shed some light on their appearance.

We begin by studying Dedekind--Carlitz sums through the generating function attached to the lattice points of a rational polyhedron $\P \subseteq \R^d$, namely, the \emph{integer-point transform}
\[
  \sigma_\P \left( \z \right) = \sigma_\P \left( z_1, z_2, \dots z_d \right) := \sum_{ \m \in {\P \cap \Z^d} } \z^\m \, ,
\]
where $\z^\m = z_1^{ m_1 } z_2^{ m_2 } \cdots z_d^{ m_d }$.
Our goals in this paper are as follows:
\begin{itemize}
  \item We show that Dedekind--Carlitz polynomials appear naturally in integer-point transforms of rational cones (Section~\ref{carlitzappearssection}).
  \item We give novel \emph{geometric} proofs of Theorem~\ref{carlitzreciprocity}, some of its generalizations, and some new reciprocity theorems (Sections~\ref{carlitzrecsection} and~\ref{variationssection}).
  \item We show that our geometric setup immediately implies that (higher-dimensional) Dedekind--Carlitz polynomials can be computed in polynomial time (Section~\ref{complexitysection}).
  \item We realize the equivalence of Dedekind--Carlitz polynomials and the integer-point transform of a two-dimensional analogue of the Mordell--Pommersheim tetrahedron (Section~\ref{trianglesection}).
  \item We give an intrinsic geometric reason why Dedekind sums appear in Theorem~\ref{thm.L} by applying Brion's decomposition theorem \cite{brion} to the Mordell--Pommersheim tetrahedron (Section~\ref{mordellsection}).
\end{itemize}

While we believe that this paper constitutes the first fundamental study of Carlitz--Dedekind sum through a geometric setup, a part of this setup can implicitly be found in papers by Solomon \cite{solomondedsum} and Chapman \cite{chapmandedekind} on generalized Dedekind sums, as well as in \cite{hardysums}.


\section{Polyhedral Cones Give Rise to Dedekind--Carlitz Polynomials}\label{carlitzappearssection}

We start by decomposing the first quadrant $\R_{ \ge 0 }^2$ into two cones, namely
\begin{align*}
  \K_1 &= \left\{ \lambda_1 (0,1) + \lambda_2 (a,b) : \, \lambda_1, \lambda_2 \ge 0 \right\} , \\
  \K_2 &= \left\{ \lambda_1 (1,0) + \lambda_2 (a,b) : \, \lambda_1 > 0, \lambda_2 \ge 0 \right\} .
\end{align*}
Thus $\K_1$ is closed, $\K_2$ is half-open, and $\R_{ \ge 0 }^2$ is the disjoint union of $\K_1$ and $\K_2$.
Let's compute their integer-point transforms. By a simple tiling argument (see, for example, \cite[Chapter~3]{ccd}),
\[
  \sigma_{ \K_1 } (u,v)
  = \sigma_{ \Pi_1 } (u,v) \left( \sum_{ j \ge 0 } v^j \right) \left( \sum_{ k \ge 0 } u^{ k a } v^{ k b } \right)
  = \frac{ \sigma_{ \Pi_1 } (u,v) }{ \left( 1 - v \right) \left( 1 - u^a v^b \right) } \, ,
\]
where 
\[
  \Pi_1 = \left\{ \lambda_1 (0,1) + \lambda_2 (a,b) : \, 0 \le \lambda_1, \lambda_2 < 1 \right\}
\]
is the \emph{fundamental parallelogram} of the cone $\K_1$.
Analogously, we can write
\[
  \sigma_{ \K_2 } (u,v)
  = \frac{ \sigma_{ \Pi_2 } (u,v) }{ \left( 1 - u \right) \left( 1 - u^a v^b \right) } \, ,
\]
where 
\[
  \Pi_2 = \left\{ \lambda_1 (1,0) + \lambda_2 (a,b) : \, 0 < \lambda_1 \le 1 , \, 0 \le \lambda_2 < 1 \right\} .
\]
Note that we need to include different sides of the half-open parallelograms $\Pi_1$ and $\Pi_2$.

\begin{figure}[ht] 
  \begin{center}
    \includegraphics[totalheight=2in]{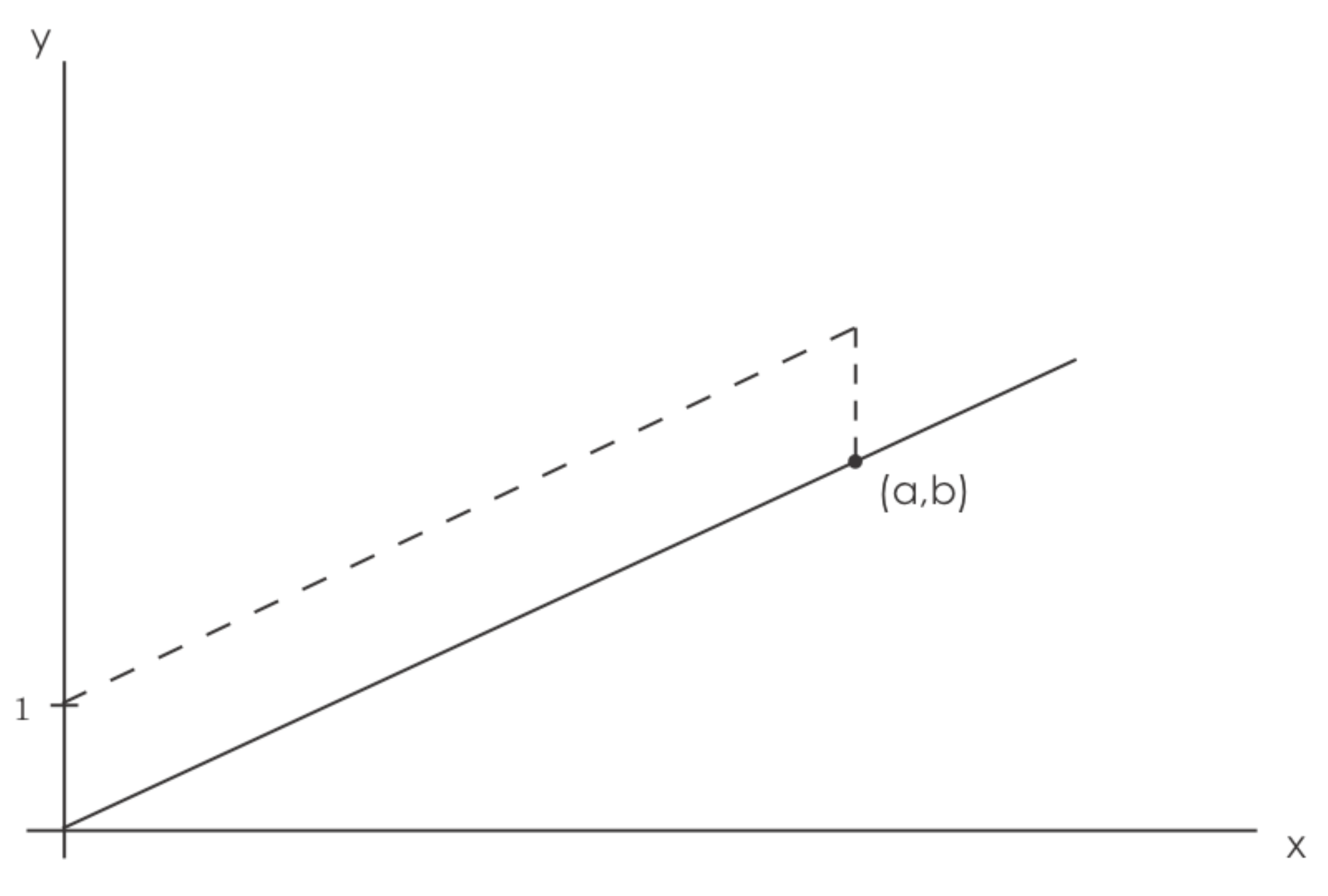}
  \end{center}
  \caption{The fundamental parallelogram $\Pi_1$.}\label{fundparal}
\end{figure}

Now we list the integer points in the half-open parallelepiped $\Pi_1$. We may assume that $a$ and $b$ are relatively prime. Then, since $\Pi_1$ has height 1,
\[
  \Pi_1 \cap \Z^2 = \left\{ (0,0), \left( k, \fl{ \frac{kb}{a} } +1 \right) : \, 1 \leq k \leq a-1, \, k \in \Z \right\} ,
\]
whence
\[
  \sigma_{K_1} (u,v) = \frac{ 1 + \sum_{ k=1 }^{ a-1 } u^k v^{ \fl{ \frac{kb}{a} } +1 } }{ (1 - v)(1 - u^a v^b) } =  \frac{ 1+ u v \, \car \left( v,u; b,a \right) }{ (v-1) \left( u^a v^b - 1 \right) } \, .
\]

We obtain the integer-point transform for $K_2$ in the same way, carefully adjusting our sums for the half-open parallelepiped $\Pi_2$:
\[
  \sigma_{K_2} (u,v) = \frac{ u+ \sum_{ k=1 }^{ b-1 } v^k u^{ \fl{ \frac{ka}{b} } + 1 } }{ (1 - u)(1 - u^a v^b) } = \frac{ u + u v \, \car \left(u,v; a,b \right) }{ (u-1) \left( u^a v^b - 1 \right) } \, .
\]


\section{Carlitz Reciprocity}\label{carlitzrecsection}

The reciprocity theorem for Dedekind--Carlitz sums now follows almost instantly from our geometric setup.

\begin{proof}[Proof of Theorem \ref{carlitzreciprocity}]
We have constructed two cones $\K_1$ and $\K_2$ such that $\R_{ \ge 0 }^2 = \K_1 \sqcup \K_2$ as a disjoint union. In the language of integer-point transforms, this means
\begin{equation} \label{2dimrec}
  \sigma_{ \R_{ \ge 0 }^2 } (u,v) = \sigma_{ \K_1 } (u,v) + \sigma_{ \K_2 } (u,v) \, .
\end{equation}
We just computed the rational generating functions on the right-hand side, and the integer-point transform for the first quadrant is simple:
$\sigma_{ \R_{ \ge 0 }^2 } (u,v) = \frac{1}{ (1 - u)(1 - v) }$.
Thus \eqref{2dimrec} becomes
\[
  \frac{1}{ (u-1)(v-1) } = \frac{ 1+ u v \, \car \left( v,u; b,a \right) }{ (v-1) \left( u^a v^b - 1 \right) } + \frac{ u+ u v \, \car \left(u,v; a,b \right) }{ (u-1) \left( u^a v^b - 1 \right) } \, ,
\]
which yields the identity of Theorem \ref{carlitzreciprocity} after clearing denominators.
\end{proof}

Our new, geometric proof of Carlitz's reciprocity theorem has a natural generalization to higher dimensions, which yields the reciprocity identity for the \emph{higher-dimensional Dedekind--Carlitz polynomials}
\[ 
  \car \left(u_1, u_2, \dots, u_n ; a_1, a_2, \dots, a_n \right) := \sum_{k=1}^{a_n -1} u_1^{ \fl{ \frac{k a_1}{a_n} } } u_2^{ \fl{ \frac{k a_2}{a_n} } } \cdots u_{n-1}^{ \fl{ \frac{k a_{n-1}}{a_n} } } u_n^{k-1} ,
\]    
where $u_1, u_2, \dots, u_n$ are indeterminates and $a_1, a_2, \dots, a_n$ are positive integers.
The reciprocity theorem for these polynomials is due to Berndt and Dieter \cite{berndtdieter}, and we give a novel proof using essentially the same geometric picture as in our previous proof.

\begin{theorem}[Berndt--Dieter]\label{carlitz reciprocity n-dim}
If $a_1, a_2, \dots, a_n$ are pairwise relatively prime positive integers, then  
\begin{align*}
  &\left(u_n -1\right) \car \left(u_1, u_2, \dots, u_n ; a_1, a_2, \dots, a_n \right) \\
  &\qquad +  \left(u_{n-1} -1\right) \car \left(u_n, u_1, \dots, u_{n-2}, u_{n-1} ; a_n, a_1, \dots, a_{n-2}, a_{n-1} \right) \\
  &\qquad + \cdots + \left(u_1 -1\right) \car \left(u_2, u_3, \dots, u_n, u_1 ; a_2, a_3, \dots, a_n, a_1 \right) \\
  &\qquad = u_1^{a_1 -1} u_2^{a_2-1} \cdots u_n^{a_n-1} - 1 \, .
\end{align*}
\end{theorem}

\begin{proof}
Analogous to our proof of Theorem \ref{carlitzreciprocity}, we construct a single ray in $n$-dimensional space, and then we decompose the non-negative orthant into $n$ cones as follows. Let $\a:=(a_1, a_2, \dots, a_n) \in \R^n$, denote the $j\th$ unit vector by $\e_j$, and define
\begin{align*}
  \K_1 &= \left \{ \lambda_2 {\e}_2 + \lambda_3 {\e}_3 + \cdots + \lambda_n {\e}_n + \lambda \a : \, \lambda_2, \dots, \lambda_n, \lambda \ge 0 \right \}, \\
  \K_2 &= \left \{ \lambda_1 {\e}_1 + \lambda_3 {\e}_3 + \cdots + \lambda_n {\e}_n + \lambda \a : \, \lambda_1 > 0, \, \lambda_3, \dots, \lambda_n, \lambda \ge 0 \right \}, \\
  &\vdots \\
  \K_j &= \left\{  \begin{array}{l}
    \lambda_1 {\e}_1 + \cdots + \lambda_{j-1} {\e}_{j-1} +\lambda_{j+1} {\e}_{j+1} + \cdots + \lambda_n {\e}_n + \lambda \a : \\
    \lambda_1, \dots, \lambda_{ j-1 } > 0 , \, \lambda_{ j+1 } , \dots, \lambda_n, \lambda \ge 0
  \end{array}  \right\}, \\
  &\vdots \\
  \K_n &= \left \{ \lambda_1 {\e}_1 + \lambda_2 {\e}_2 + \cdots + \lambda_{n-1} {\e}_{n-1} + \lambda \a : \, \lambda_1, \dots, \lambda_{ n-1 } > 0, \, \lambda \ge 0 \right \} .
\end{align*}
The fundamental parallelepiped of $\K_j$ is
\[
  \Pi_j = \left\{  \begin{array}{l}
    \lambda_1 {\e}_1 + \cdots + \lambda_{j-1} {\e}_{j-1} +\lambda_{j+1} {\e}_{j+1} + \cdots + \lambda_n {\e}_n + \lambda \a : \\
    0 < \lambda_1, \dots, \lambda_{ j-1 } \le 1 , \, 0 \le \lambda_{ j+1 } , \dots, \lambda_n, \lambda < 1
  \end{array}  \right\} .
\]
Thus a point in $\Pi_j$ will look like
\[
  \left( \lambda_1 + \lambda a_1, \dots, \lambda_{ j-1 } + \lambda a_{ j-1 } , \lambda a_j, \lambda_{ j+1 } + \lambda a_{ j+1 } , \dots, \lambda_n + \lambda a_n \right) ,
\]
and a slice of this parallelepiped at $x_j = k$ (where $1 \le k \le a_j-1$) will contain
\[
  \left( \fl{ \frac{ka_1}{a_j}} + 1, \dots, \fl{ \frac{ka_{j-1}}{a_j}} + 1, k, \fl{ \frac{ka_{j+1}}{a_j}} + 1, \dots, \fl{ \frac{ka_n}{a_j}} + 1 \right)
\]
as the only integer point.
Note also that the integer point $\e_1 + \dots + \e_{ j-1 }$ is in $\Pi_j$.
Combining this information as in Section \ref{carlitzappearssection} yields the integer-point transform
\begin{align*}
  &\sigma_{\K_j} (\u) = \frac{ \sigma_{ \Pi_j } (\u) }{ \left( 1-\u^\a \right) \left( 1-u_1 \right) \cdots \left( 1-u_{ j-1 } \right) \left( 1-u_{ j+1 } \right) \cdots \left( 1-u_n \right) } \\
  &\qquad = \frac{ u_1 u_2 \cdots u_{j-1} + \sum_{k=1}^{a_j -1} u_1^{\fl{ \frac{ka_1}{a_j}} + 1} \cdots u_{j-1}^{\fl{ \frac{ka_{j-1}}{a_j}} + 1} u_j^k u_{j+1}^{\fl{ \frac{ka_{j+1}}{a_j}} + 1} \cdots u_n^{\fl{ \frac{ka_n}{a_j}} + 1} }{ \left( 1-\u^\a \right) \left( 1-u_1 \right) \cdots \left( 1-u_{ j-1 } \right) \left( 1-u_{ j+1 } \right) \cdots \left( 1-u_n \right) } \\
  &\qquad = \frac{ u_1 u_2 \cdots u_{j-1} + u_1 u_2 \cdots u_n \, \car \left( u_1, \dots, u_{ j-1 }, u_{ j+1 } , \dots, u_n, u_j ; a_1, \dots, a_{ j-1 }, a_{ j+1 } , \dots, a_n, a_j \right) }{ \left( 1-\u^\a \right) \left( 1-u_1 \right) \cdots \left( 1-u_{ j-1 } \right) \left( 1-u_{ j+1 } \right) \cdots \left( 1-u_n \right) } \, ,
\end{align*}
where $\u := \left( u_1, u_2, \cdots, u_n \right)$ and $\u^{\a} := u_1^{a_1} u_2^{a_2} \cdots u_n^{a_n}$. 
Since $\bigcup_{j=1}^n \K_j = \R_{ \ge 0 }^n$ is a disjoint union of the nonnegative $n$-dimensional orthant,
\[
  \sigma_{\K_1} (\u) + \sigma_{\K_2} (\u) + \cdots + \sigma_{\K_n} (\u) = \sigma_{\R_{ \ge 0 }^n} (\u) = \frac{ 1 }{ \left( 1-u_1 \right) \cdots \left( 1-u_n \right) } \, .
\]
Theorem \ref{carlitz reciprocity n-dim} follows upon clearing denominators in this identity.
\end{proof}

To virtually every theorem in this paper, there exist translate companions, i.e., we can shift the cones involved in our proofs by a fixed vector. This gives rise to shifts in the greatest-integer functions, and the resulting Carlitz sums are polynomial analogues of \emph{Dedekind--Rademacher sums} \cite{rademacherdedekind}. For the sake of clarity of exposition, we only give integer-vertex versions of our cones, but the reader should keep in mind that arbitrary vertices do, in principal, not cause any additional problems.


\section{Computational Complexity}\label{complexitysection}

Dedekind's reciprocity law \eqref{dedreci} together with the identity $\ded (a, b) = \ded (c, b)$, if $a \equiv c \bmod b$, allows us to compute the classical Dedekind sum in a Euclidean-algorithm style and thus very efficiently, namely in linear time. (Computational complexity is measured in terms of the input length, e.g., in this case $\log a + \log b$.)
We do not know how to apply a similar reasoning to the Dedekind--Carlitz polynomials via Carlitz's reciprocity law Theorem \ref{carlitzreciprocity}; however, the following central theorem of Barvinok \cite{barvinokalgorithm} allows us to deduce that (higher-dimensional) Dedekind--Carlitz polynomials can be computed efficiently.

\begin{theorem}[Barvinok]
In fixed dimension, the integer-point transform $\sigma_\P \left( z_1, z_2, \dots, z_d \right)$ of a rational polyhedron $\P$ can be computed as a sum of rational functions in $z_1, z_2, \dots, z_d$ in time polynomial in the input size of $\P$.
\end{theorem}

Note that Barvinok's theorem says, in particular, that the rational functions whose sum represents $\sigma_\P (\z)$ are \emph{short}, i.e., the set of data needed to output this sum of rational functions is of size that is polynomial in the input size of $\P$.
The application of Barvinok's theorem to any of the cones appearing in the proof of Theorem~\ref{carlitz reciprocity n-dim} immediately yields the following novel complexity result.

\begin{theorem}\label{carcomplexity}
For fixed $n$, the higher-dimensional Dedekind--Carlitz polynomial \\
$\car \left(u_1, u_2, \dots, u_n ; a_1, a_2, \dots, a_n \right)$ can be computed in time polynomial in the size of $a_1, a_2, \dots, a_n$.
\end{theorem}

In particular, this result says that there is a more economical way to write the ``long" polynomial $\car \left(u_1, u_2, \dots, u_n ; a_1, a_2, \dots, a_n \right)$ as a short sum of rational functions.
Theorem~\ref{carcomplexity} implies that any Dedekind-like sum that can be efficiently derived from Dedekind--Carlitz polynomials (e.g., by applying differential operators) can also be computed efficiently.


\section{Variations on a Theme}\label{variationssection}

The geometry we have used so far is very simple: it is based on one ray in space.
Naturally, this can be extended in numerous ways. To illustrate one of them, consider two rays through the points $(a,b)$ and $(c,d)$ in the first quadrant of the plane, where $a,b,c,d$ are positive integers.  This construction decomposes the first quadrant into three rational cones as shown in Figure~\ref{tworays}.

\begin{figure}[ht] 
  \begin{center}
    \includegraphics[totalheight=2in]{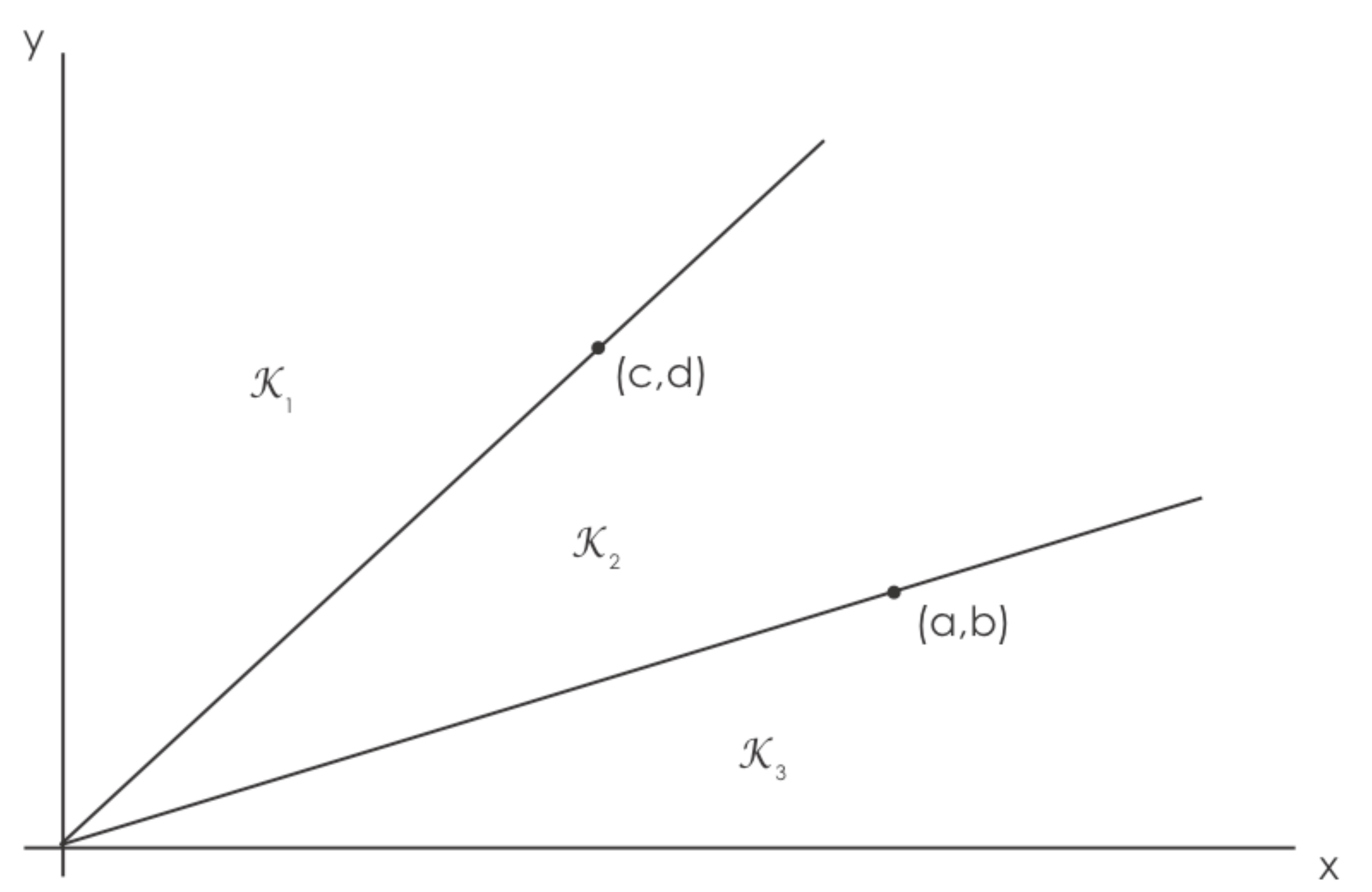}
  \end{center}
  \caption{Two-ray decomposition of the first quadrant.}\label{tworays.pdf}
\end{figure}

The integer-point transform of each of the two exterior cones, $\K_1$ and $\K_3$, is easily computed in the same manner as in Section \ref{carlitzappearssection}.  However, the cone in the middle, $\K_2$, is bounded on either side by a non-unit vector, and therefore the interior points of the fundamental parallelogram of this cone are not trivial to list. We will use unimodular transformations to compute $\sigma_{ \K_2 }$.
This construction leads to a new reciprocity identity (Theorem~\ref{carlitz abcd} below). We start by computing the integer-point transform for a general rational cone in $\R^2$, i.e., $\K_2$ in Figure \ref{tworays}.

\begin{lemma}\label{general2dimcone}
Suppose $a,b,c,d$ are positive integers such that $ad > bc$ and
$\gcd(a,b) = \gcd(c,d) = 1$, and let $x, y \in \Z$ such that $ax + by
= 1$. Then the cone 
$\K := \left\{ \lambda (a,b) + \mu (c,d) : \, \lambda, \mu \geq 0 \right\}$
has the integer-point transform
\[
  \sigma_{ \K } (u,v)
  = \frac{ 1 + u^{ a-y } v^{ b+x } \, \car \left( u^a v^b , u^{-y} v^x ; cx+dy , ad-bc \right) }{ \left( u^a v^b - 1 \right) \left( u^c v^d - 1 \right) } \, .
\]
\end{lemma}

\begin{proof}
To compute the generating function $\sigma_{\K}$, note that the linear transformation given by
\[
  M := \left( \begin{array}{cc} 
    x & y \\
    -b & a
  \end{array} \right) 
\]
maps $\K$ to the cone $M(\K)$ generated by $(1,0)$ and $(cx+dy, ad-bc)$, whose integer-point transform we know how to compute from Section~\ref{carlitzappearssection}.
Since
\[
  M^{ -1 } = \left( \begin{array}{cc} 
    a & -y \\
    b & x
  \end{array} \right) 
\]
we obtain
\begin{align*}
  \sigma_{ \K } (u,v)
  &= \sum_{ (m,n) \in \K \cap \Z^2 } u^m v^n
  = \sum_{ (m,n) \in M(\K) \cap \Z^2 } u^{ am-yn } v^{ bm+xn }
  = \sum_{ (m,n) \in M(\K) \cap \Z^2 } \left( u^a v^b \right)^m \left( u^{-y} v^x \right)^n \\
  &= \frac{ 1 + u^{ a-y } v^{ b+x } \car \left( u^a v^b , u^{-y} v^x ; cx+dy , ad-bc \right) }{ \left( u^a v^b - 1 \right) \left( u^c v^d - 1 \right) } \, .  \qedhere
\end{align*}
\end{proof}

\begin{theorem}\label{carlitz abcd}
Suppose $a,b,c,d$ are positive integers such that $ad > bc$ and $\gcd(a,b) = \gcd(c,d) = 1$, and let $x, y \in \Z$ such that $ax + by = 1$. Then
\begin{align*}
  &uv(u-1) \left( u^a v^b -1 \right) \car \left( v, u; d, c \right) 
  + uv (v-1) \left( u^c v^d - 1 \right) \car \left( u, v; a, b \right) \\
  &\qquad + u^{ a-y } v^{ b+x } \left(  u-1 \right) \left( v-1 \right) \car \left( u^a v^b , u^{-y} v^x ; cx+dy , ad-bc \right) \\
  &\qquad = u^{ a+c } v^{ b+d } - u^a v^b (uv-v+1) - u^c v^d (uv-u+1) + uv \, .
\end{align*}
\end{theorem}

\begin{proof}
Again, we prove the theorem geometrically.  Let
\begin{align*}
  &\K_1 = \left\{ \lambda_2 {\e}_2 + \lambda_{cd} (c,d) : \, \lambda_2 > 0, \, \lambda_{cd} \geq 0 \right\} , \\
  &\K_2 = \left\{ \lambda_{ab} (a,b) + \lambda_{cd} (c,d) : \, \lambda_{ab}, \lambda_{cd} \geq 0 \right\} , \\
  &\K_3 = \left\{ \lambda_1 {\e}_1 + \lambda_{ab} (a,b) : \, \lambda_1 > 0, \, \lambda_{ab} \geq 0 \right\} ,
\end{align*}
so that $\K_2$ is closed and $\K_1$ and $\K_3$ are half-open, and $\K_1 \sqcup \K_2 \sqcup \K_3 = \R_{ \ge 0 }^2$ is a disjoint union of the first quadrant.  With the methods introduced in Section \ref{carlitzappearssection}, the integer-point transforms of $\K_1$ and $\K_3$ are
\begin{align*}
  &\sigma_{\K_1} (u,v) = \frac{ v + \sum_{k=1}^{c-1} u^k v^{ \fl{\frac{kd}{c}}+1} }{ (1-v) \left( 1-u^c v^d \right) } = \frac{ v + u v \, \car \left(v,u; d,c \right) }{ (v-1) \left( u^c v^d -1 \right) } \, , \\
  &\sigma_{\K_3} (u,v) = \frac{ u + \sum_{k=1}^{b-1} v^k u^{ \fl{\frac{ka}{b}}+1} }{ (1-u) \left( 1-u^a v^b \right) } = \frac{ u + u v \, \car \left( u,v ; a,b \right) }{ (u-1) \left( u^a v^b -1 \right) } \, ,
\end{align*}
and the integer-point transform of $\K_2$ was computed in Lemma~\ref{general2dimcone}.
By our construction,
\[
  \sigma_{\K_1} (u,v) + \sigma_{\K_2} (u,v) + \sigma_{\K_3} (u,v) = \sigma_{\R_{ \ge 0 }^2} (u,v) = \frac{ 1 }{ (1-u)(1-v) }  \, ,
\]
from which Theorem \ref{carlitz abcd} follows upon clearing denominators.
\end{proof}

Theorem~\ref{carlitz abcd} is the polynomial analogue of the following result due to Pommersheim \cite[Theorem~7]{pommersheim}:

\begin{theorem}[Pommersheim]\label{pomm3term}
Suppose $a,b,c,d$ are positive integers such that $\gcd(a,b) = \gcd(c,d) = 1$, and let $x, y \in \Z$ such that $ax + by = 1$. Then
\begin{align*}
  \ded (a,b) + \ded(c,d) = \ded (cx-dy, ad+bc) - \frac 1 4 + \frac{ 1 }{ 12 } \left( \frac{ b }{ d (ad+bc) } + \frac{ d }{ b (ad+bc) } + \frac{ ad+bc }{ bd } \right) .
\end{align*}
\end{theorem}

This, in turn, generalizes Rademacher's three-term reciprocity theorem \cite{rademacherthreeterm}:

\begin{corollary}[Rademacher]\label{radcor}
If $a, b, c$ are pairwise coprime positive integers then
\[
  \ded \left( a b^{ -1 } , c \right) + \ded \left( c a^{ -1 } , b \right) + \ded \left( b c^{ -1 } , a \right) = - \frac 1 4 + \frac{ 1 }{ 12 } \left( \frac{ a }{ bc } + \frac{ b }{ ca } + \frac{ c }{ ab } \right) ,
\]
where $a^{ -1 } a \equiv 1 \bmod b$, $b^{ -1 } b \equiv 1 \bmod c$, and $c^{ -1 } c \equiv 1 \bmod a$.
\end{corollary}

The case $n=3$ of Berndt--Dieter's Theorem~\ref{carlitz reciprocity n-dim} is the polynomial analogue of this three-term law. 
Clearly, Dedekind's reciprocity theorem \eqref{dedreci} is implied by Corollary~\ref{radcor}, and Girstmair showed \cite{girstmair3term} that \eqref{dedreci} implies Theorem~\ref{pomm3term}, so that, in fact, Dedekind's, Rademacher's, and Pommersheim's reciprocity theorems are equivalent. 
Such an equivalence can clearly not hold on the polynomial level, but one could ask whether Theorems \ref{carlitzreciprocity} and \ref{carlitz abcd} are equivalent.

Theorem~\ref{carlitz abcd} could be generalized in several ways, e.g., to higher dimensions or to more than three cones in dimension 2 (this yields a Carlitz polynomial analogue of \cite[Theorem~8]{pommersheim}), but we digress.

Theorem \ref{carlitz abcd} simplifies when $ad - bc = 1$: then the third Dedekind--Carlitz polynomial disappears. Geometrically, this stems from the fact that the cone $\K_2$ is \emph{unimodular}, i.e., its fundamental parallelogram contains only the origin. 

\begin{corollary}
If $a,b,c,d$ are positive integers such that $ad - bc = 1$, then
\begin{align*}
  &(u-1) \left( u^a v^b -1 \right) \car \left( v, u; d,c \right) + (v-1) \left( u^c v^d -1 \right) \car \left( u,v; a,b \right) \\
  &\qquad = u^{a+c-1} v^{b+d-1} - u^a v^b -u^c v^d + u^{a-1} v^b + u^c v^{d-1}  - u^{a-1} v^{b-1} - u^{c-1} v^{d-1} +1 \, .
\end{align*}
\end{corollary}

This is the polynomial analogue of the following result due to Rademacher~\cite{rademachersummen}.

\begin{corollary}[Rademacher]\label{rademacher abcd}
If $a,b,c,d$ are positive integers such that $ad-bc=1$, then
\[
  \ded (a,b) + \ded (d,c) = - \frac{1}{2} + \frac{1}{12} \left( \frac{a}{b} + \frac{a}{c} + \frac{d}{b} + \frac{d}{c} \right) .
\]
\end{corollary}


\section{Dedekind--Carlitz Polynomials as Integer-Point Transforms of Triangles}\label{trianglesection}

We have shown that Dedekind--Carlitz polynomials appear as natural ingredients of integer-point transforms of cones. In this light, it should come as no surprise that the following conic decomposition theorem of Brion \cite{brion} will prove useful. 
We remind the reader that a \emph{rational polytope} is the convex hull of finitely many rational points in $\R^d$.
If $\v$ is a vertex of the polytope $\P$, then the \emph{vertex cone} $\K_\v$ is the smallest cone with apex $\v$ that contains $\P$.

\begin{theorem}[Brion] \label{brionthm}
Suppose $\P$ is a rational convex polytope. Then we have the following identity of rational functions: 
\[
  \sigma_\P (\z) = \sum_{ \v } \sigma_{\K_\v} (\z) \, ,
\]
where the sum is over all vertices of $P$.
\end{theorem}  

Brion's theorem allows us to give a novel expression for the Dedekind--Carlitz polynomial as the integer-point transform of a certain triangle.

\begin{theorem}\label{thm.brion2d}
Let $a$ and $b$ be relatively prime positive integers and $u$ and $v$ be indeterminates.  If $\Delta$ is the triangle with vertices $(0,0)$, $(a,0)$, and $(0,b)$, then the Dedekind--Carlitz polynomial $\car \left( \frac{1}{u}, v; a, b \right)$ and the integer-point transform of $\Delta$ are related in the following manner:
\[
  (u-1) \, \sigma_\Delta (u,v) = u^a v \, \car \left( \frac{1}{u}, v; a, b \right) + u \left( u^a + v^b \right) - \frac{ v^{b+1} - 1 }{ (v-1) } \, .
\]
\end{theorem}

\begin{proof}
The triangle $\Delta$ comes with the vertex cones
\begin{align*}
  \K_1 &= \left \{ j (1,0) + k (0,1) : \, j,k \geq 0 \right \}, \\
  \K_2 &= \left \{ (a,0) + j (-1,0) + k (-a,b) : \, j,k \geq 0 \right \}, \\  
  \K_3 &= \left \{ (0,b) + j (a,-b) + k (0,-1) : \, j,k \geq 0 \right \}.
\end{align*}
Once more we apply the methods of Section \ref{carlitzappearssection} and obtain the integer-point transforms
\begin{align*}
  \sigma_{\K_1} (u,v) &= \frac{1}{ (1-u)(1-v) } \, , \\
  \sigma_{\K_2} (u,v) &= u^a \frac{1 + \sum_{k=1}^{b-1} v^k u^{ \fl{ - \frac{ka}{b} } } }{ \left( 1 - u^{ -1 }  \right) \left( 1- u^{-a} v^b \right) } = - \frac{ u^{a+1} + u^a v \, \car \left( \frac{1}{u}, v; a, b \right) }{ (u-1) \left( u^{-a} v^b -1 \right) } \, , \\
  \sigma_{\K_3} (u,v) &= v^b \frac{1 + \sum_{k=1}^{a-1} u^k v^{ \fl{ - \frac{kb}{a} } } }{ \left( 1 - v^{-1} \right) \left( 1- u^a v^{-b} \right) } = - \frac{ v^{b+1} + u v^b \, \car \left( \frac{1}{v}, u; b, a \right) }{ (v-1) \left( u^a v^{-b} -1 \right) } \, .
\end{align*}
Now Brion's Theorem \ref{brionthm} gives
\[
  \sigma_\Delta (u,v) =  \frac{1}{ (u-1)(v-1) } - \frac{ u^{a+1} + u^a v \, \car \left( \frac{1}{u}, v; a, b \right) }{ (u-1) \left( u^{-a} v^b -1 \right) } - \frac{ v^{b+1} + u v^b \, \car \left( \frac{1}{v}, u; b, a \right) }{ (v-1) \left( u^a v^{-b} -1 \right) }
\]
or, after some manipulations,
\begin{multline}
  (u-1)(v-1) \left( u^a - v^b \right) \sigma_\Delta (u,v) = u^a - v^b + u^{2a+1} (v-1) - v^{2b+1} (u-1) \\
  \qquad + u^{2a} v (v-1) \, \car \left( \frac{1}{u}, v; a, b \right) - u v^{2b} (u-1) \, \car \left( \frac{1}{v}, u; b, a \right) \, . \label{brionrewrite}
\end{multline}
The two Dedekind--Carlitz polynomials in this expression are actually related by reciprocity: after a simple change of variables, Carlitz's Theorem~\ref{carlitzreciprocity} implies
\[
  u v^b (u-1) \, \car \left( \frac{1}{v}, u; b, a \right) - u^a v (v-1) \, \car \left( \frac{1}{u}, v; a, b \right) = u^a v - u v^b ,
\]
and substituting this into \eqref{brionrewrite} finishes the proof.
\end{proof}

We remark that we could also 
deal with triangles with rational vertices. The setup is exactly the same, however, the formulas become quite a bit messier.


\section{The Mordell--Pommersheim Tetrahedron}\label{mordellsection}

In this final section we derive a generating-function equivalent to Theorem~\ref{thm.L}.
%
%
It is natural to extend the application of Brion's theorem in the last section to polytopes in higher dimensions, which we will do with the tetrahedron $\T$. Since our ultimate goal is a geometric proof of Theorem~\ref{thm.L}, we will apply Brion's theorem to the dilate $t \T$ for any positive integral $t$. As we will see, the following higher-dimensional extensions of Dedekind--Carlitz polynomials will appear naturally in the integer-point transforms of the vertex cones of $\T$.
For positive integers $a,b,c$, and indeterminates $u,v,w$, we define the \emph{Dedekind--Rademacher--Carlitz (DRC) sum}
\[ 
  \drc (u,v,w;a,b,c) := \sum_{k=0}^{c-1} \sum_{j=0}^{b-1} u^{ \fl{ \frac{ja}{b} + \frac{ka}{c} } } v^j w^k.
\]
The DRC sums are the main ingredients for the integer-point transform of $\T$.

\begin{theorem}\label{thm.brion3d}
Let $\T$ be the convex hull of $(a,0,0)$, $(0,b,0)$, $(0,0,c)$, and $(0,0,0)$, where $a, b$ and $c$ are pairwise relatively prime positive integers. Then
\begin{align*}
  &(u-1)(v-1)(w-1) \left( u^a - v^b \right) \left( u^a - w^c \right) \left( v^b - w^c \right) \sigma_{t \T} (u,v,w) \\
  &\quad = u^{(t+2)a} (v-1)(w-1) \left( v^b - w^c \right) \left( (u-1) + \drc \left( u^{-1},v,w;a,b,c \right) \right) \\
  &\qquad - v^{(t+2)b} (u-1)(w-1) \left( u^a - w^c \right) \left( (v-1) + \drc \left( v^{-1},u,w;b,a,c \right) \right) \\
  &\qquad + w^{(t+2)c} (u-1)(w-1) \left( u^a - v^b \right) \left( (w-1) + \drc \left( w^{-1},u,v;c,a,b \right) \right) \\
  &\qquad - \left( u^a - v^b \right) \left( u^a - w^c \right) \left( v^b - w^c \right).
\end{align*}
\end{theorem}  

\begin{proof}
The tetrahedron $t\T$ has the vertex cones
\begin{align*}
  \K_0 &= \left \{ \lambda_1 (1,0,0) + \lambda_2 (0,1,0) + \lambda_3 (0,0,1) : \, \lambda_1, \lambda_2, \lambda_3 \geq 0 \right \}, \\
  \K_1 &= \left \{ (ta,0,0) + \lambda_1 (-1,0,0) + \lambda_2 (-a,b,0) + \lambda_3 (-a,0,c) : \, \lambda_1, \lambda_2, \lambda_3 \geq 0 \right \}, \\
  \K_2 &= \left \{ (0,tb,0) + \lambda_1 (a,-b,0) + \lambda_2 (0,-1,0) + \lambda_3 (0,-b,c) : \, \lambda_1, \lambda_2, \lambda_3 \geq 0 \right \}, \\
  \K_3 &= \left \{ (0,0,tc) + \lambda_1 (a,0,-c) + \lambda_2 (0,b,-c) + \lambda_3 (0,0,-1) : \, \lambda_1, \lambda_2, \lambda_3 \geq 0 \right \} ,
\end{align*}  
and we will ultimately apply Brion's Theorem~\ref{brionthm} to obtain
\begin{equation}\label{3dimbriontet}
  \sigma_{ t \T } (\z) = \sigma_{ \K_0 } (\z) + \sigma_{ \K_1 } (\z) + \sigma_{ \K_2 } (\z) + \sigma_{ \K_3 } (\z) \, ,
\end{equation}
which will yield Theorem~\ref{thm.brion3d}.

The computation of the integer-point transforms of the vertex cones should be routine by now; we will derive only one of them in detail.
The vertex cone $\K_1$ has the fundamental parallelepiped
\[
  \Pi = (ta,0,0) + \left \{ \lambda_1 (-1,0,0) + \lambda_2 (-a,b,0) + \lambda_3
    (-a,0,c) : \, 0 \leq \lambda_1, \lambda_2, \lambda_3 < 1 \right \} .
\]
There exists exactly one integer point in $\Pi$ at each integer $j$ along the $y$-axis and $k$ along the $z$-axis.  In other words, for $j,k \in \Z_{\geq 0}$, $(x,j,k) \in \Pi$ is written as $(x,j,k) = \left( - \lambda_1 - \lambda_2 a - \lambda_3 a, \lambda_2 b, \lambda_3 c \right)$ and hence $\lambda_2 = \frac{j}{b}$ and $\lambda_3 = \frac{k}{c}$. Thus
\[
  \Pi \cap \Z^3 = \left\{ \left( \fl{ - \frac{ja}{b} - \frac{ka}{c} }, j, k \right) : \, j,k = 0,1, \dots, b-1 \right\}
\]
and
\begin{align*}
  \sigma_{\K_1} (u,v,w) &= u^{ta} \left( \frac{ \sum_{k=0}^{c-1} \sum_{j=0}^{b=1} u^{ \fl{ - \frac{ja}{b} - \frac{ka}{c} } } v^j w^k }{ \left( 1-u^{-1} \right) \left( 1-u^{-a} v^b \right) \left( 1-u^{-a} w^c \right) } \right) \\
  &= u^{ta} \left( \frac{ 1 + \sum_{k=0}^{c-1} \sum_{j=0}^{b=1} u^{ - \fl{ \frac{ja}{b} + \frac{ka}{c} } - 1 } v^j w^k - u^{-1} }{ u^{-2a-1} \left( u-1 \right) \left( u^a - v^b \right) \left( u^a - w^c \right) } \right) \\
  &= \frac{ u^{(t+2)a} \left[ (u-1) + \drc \left( u^{-1}, v, w; a, b, c \right) \right] }{ \left( u-1 \right) \left( u^a - v^b \right) \left( u^a - w^c \right) } \, . 
\end{align*}
Similarly, we derive
\begin{align*}
  \sigma_{\K_2} (u,v,w) &= - \frac{ v^{(t+2)b} \left[ (v-1) + \drc \left( v^{-1}, u, w; b, a, c \right) \right] }{ \left( v-1 \right) \left( u^a - v^b \right) \left( v^b - w^c \right) } \, , \\
  \sigma_{\K_3} (u,v,w) &= \frac{ w^{(t+2)c} \left[ (w-1) + \drc \left( w^{-1}, u, v; c, a, b \right) \right] }{ \left( w-1 \right) \left( u^a - w^c \right) \left( v^b - w^c \right) } \, ,
\end{align*}
and
\[
  \sigma_{\K_0} (u,v,w) = - \frac{1}{ (u-1)(v-1)(w-1) } \, .
\]
Substituting these rational functions into \eqref{3dimbriontet} and clearing denominators gives the theorem.
\end{proof}

\begin{proof}[Proof of Theorem \ref{thm.L}]
Theorem \ref{thm.brion3d} gives $\sigma_{t\T} (u,v,w) = \frac N D$ with numerator
\begin{align*}
  N &= u^{(t+2)a} (v-1)(w-1) \left( v^b - w^c \right) \left( (u-1) + \drc \left( u^{-1},v,w;a,b,c \right) \right) \\
  &\qquad - v^{(t+2)b} (u-1)(w-1) \left( u^a - w^c \right) \left( (v-1) + \drc \left( v^{-1},u,w;b,a,c \right) \right) \\
  &\qquad + w^{(t+2)c} (u-1)(w-1) \left( u^a - v^b \right) \left( (w-1) + \drc \left( w^{-1},u,v;c,a,b \right) \right) \\
  &\qquad - \left( u^a - v^b \right) \left( u^a - w^c \right) \left( v^b - w^c \right)
\end{align*}
and denominator
$
  D = (u-1)(v-1)(w-1) \left( u^a - v^b \right) \left( u^a - w^c \right) \left( v^b - w^c \right)
$.
We obtain the Ehrhart polynomial $L_\T(t) = \sigma_{ t\T } (1,1,1)$ by taking the limit of this rational function as $u, v, w \to 1$. We need to use L'Hospital's rule and take partial derivatives with respect to $u$ once, $v$ twice, and $w$ three times in both numerator and denominator before substituting $u=v=w=1$.
This is easily done for the denominator $D$ with the result $-12 b c^2$.  

The numerator $N$ of $\sigma_{ t\T } (u,v,w)$ is not handled as easily. 
After differentiating and setting $u=v=w=1$ it becomes
\begin{align}
  &- 2 a b^2 c^3 t^3 - 6 a b^2 c^3 t^2 - 6 b c^2 t^2 - 4 a b^2 c^3 t - 18 b c^2 t \nonumber \\
  &+ 6 b c t \sum_{j = 0}^{a - 1} \sum_{k = 0}^{c - 1} \fl{ \frac{jb}{a} + \frac{kb}{c} } - 6 b c^2 t \sum_{j = 0}^{a - 1} \sum_{k = 0}^{c - 1} \fl{ \frac{jb}{a} + \frac{kb}{c} } \nonumber\\
  &- 12 b c t \sum_{j = 0}^{a - 1} \sum_{k = 0}^{c - 1} k \fl{ \frac{jb}{a} + \frac{kb}{c} } - 6 b c t \sum_{j = 0}^{a-1} \sum_{k = 0}^{b - 1} \fl{ \frac{jc}{a} + \frac{kc}{b} } \label{somewhere} \\
  &+ 24 b c^2 t \sum_{j = 0}^{a - 1} \sum_{k = 0}^{b - 1} \fl{ \frac{jc}{a} + \frac{kc}{b} } + 6 b c^2 t^2 \sum_{j = 0}^{a - 1} \sum_{k = 0}^{b - 1} \fl{ \frac{jc}{a} + \frac{kc}{b} } \nonumber \\
  &- 6 b c t \sum_{j = 0}^{a - 1} \sum_{k = 0}^{b - 1} \fl{ \frac{jc}{a} + \frac{kc}{b} } \left( 1 + \fl{ \frac{jc}{a} + \frac{kc}{b} } \right) \nonumber
\end{align}
plus numerous terms that do not depend on $t$.
Fortunately, because the constant term of $L_\T$ is $1$ (see, for example, \cite[Corollary 3.15]{ccd}), we are not concerned with these extra terms.
As before, we replace all greatest-integer functions with fractional-part functions to modify \eqref{somewhere} to
\begin{align}
  &- 2 a b^2 c^3 t^3 - 3 b^2 c^3 t^2 - 3 a b c^3 t^2 - 6 b c^2 t^2  - \frac{b^2 c^3}{a} t - 6 a b c^3 t \nonumber \\
  &- 3 b c^3 t - a c^3 t + 6 a b^2 c^2 t + 6 a b c^2 t - 18 b c^2 t - 5 a b^2 c t \nonumber \\
  &- 6 b c t \sum_{k=0}^{c - 1} \sum_{j=0}^{a-1} \left\{ \frac{jb}{a} + \frac{kb}{c} \right\} 
  + 6 b c^2 t \sum_{k=0}^{c - 1} \sum_{j=0}^{a-1} \left\{ \frac{jb}{a} + \frac{kb}{c} \right\} 
  + 12 b c t \sum_{k=0}^{c - 1} \sum_{j=0}^{a-1} k \left\{ \frac{jb}{a} + \frac{kb}{c} \right\} \label{somewhereelse} \\
  &+ 12 b c t \sum_{k=0}^{b - 1} \sum_{j=0}^{a-1} \left\{ \frac{jc}{a} + \frac{kc}{b} \right\}
  - 24 b c^2 t \sum_{k=0}^{b - 1} \sum_{j=0}^{a-1} \left\{ \frac{jc}{a} + \frac{kc}{b} \right\}  
  - 6 b c^2 t^2 \sum_{k=0}^{b - 1} \sum_{j=0}^{a-1} \left\{ \frac{jc}{a} + \frac{kc}{b} \right\} \nonumber \\ 
  &+ \frac{12 b c^2}{a} t \sum_{k=0}^{b - 1} \sum_{j=0}^{a-1} j \left\{ \frac{jc}{a} + \frac{kc}{b} \right\} 
  + 12 c^2 t \sum_{k=0}^{b - 1} \sum_{j=0}^{a-1} k\left\{ \frac{jc}{a} + \frac{kc}{b} \right\} 
  - 6 b c t \sum_{k=0}^{b - 1} \sum_{j=0}^{a-1} \left\{ \frac{jc}{a} + \frac{kc}{b} \right\}^2 . \nonumber
\end{align}
Now we use three elementary identities:
\begin{align*}
  &\sum_{k=0}^{b-1} \sum_{j=0}^{a-1} \left \{ \frac{jc}{a} + \frac{kc}{b} \right \} = \frac{ab-1}{2} \, , \\
  &\sum_{k=0}^{b-1} \sum_{j=0}^{a-1} j \left \{ \frac{jc}{a} + \frac{kc}{b} \right \} = c \, \ded (ab,c) + \frac{ab(c-1)}{4} \, , \\
  &\sum_{k=0}^{b-1} \sum_{j=0}^{a-1} \left \{ \frac{jc}{a} + \frac{kc}{b} \right \}^2 = \frac{(ab-1)(2ab-1)}{6ab} \, . 
\end{align*}
They allow us to simplify \eqref{somewhereelse} to
\begin{align*}
  &- 2 a b^2 c^3 t^3 - 3 b^2 c^3 t^2 - 3 a b c^3 t^2 - 3 b c^2 t^2 - \frac{b^2 c^3}{a} t - 3 b c^3 t - a c^3 t - 3  b^2 c^2 t \\
  &- 3 a b c^2 t - 9 b c^2 t - a b^2 c t - \frac{c}{a} t + 12 b c^2 \left( \ded (ab,c) + \ded (ac,b) + \ded (bc,a) \right) t \, .
\end{align*}
We divide by the denominator $-12 b c^2$ and add the constant term 1 to arrive at the desired formula for $L_\T(t)$.
\end{proof}



\begin{thebibliography}{10}

\bibitem{almkvist}
Gert Almkvist, \emph{Asymptotic formulas and generalized {D}edekind sums},
  Experiment. Math. \textbf{7} (1998), no.~4, 343--359.

\bibitem{barvinokalgorithm}
Alexander~I. Barvinok, \emph{A polynomial time algorithm for counting integral
  points in polyhedra when the dimension is fixed}, Math. Oper. Res.
  \textbf{19} (1994), no.~4, 769--779.

\bibitem{hardysums}
Matthias Beck, \emph{Geometric proofs of polynomial reciprocity laws of
  {C}arlitz, {B}erndt, and {D}ieter}, Diophantine analysis and related fields
  2006, Sem. Math. Sci., vol.~35, Keio Univ., Yokohama, 2006, pp.~11--18.

\bibitem{ccd}
Matthias Beck and Sinai Robins, \emph{Computing the continuous discretely:
  Integer-point enumeration in polyhedra}, Undergraduate Texts in Mathematics,
  Springer-Verlag, New York, 2007.

\bibitem{berndtdieter}
Bruce~C. Berndt and Ulrich Dieter, \emph{Sums involving the greatest integer
  function and {R}iemann-{S}tieltjes integration}, J. Reine Angew. Math.
  \textbf{337} (1982), 208--220.

\bibitem{brion}
Michel Brion, \emph{Points entiers dans les poly\`edres convexes}, Ann. Sci.
  \'Ecole Norm. Sup. (4) \textbf{21} (1988), no.~4, 653--663.

\bibitem{carlitzpolynomials}
Leonard Carlitz, \emph{Some polynomials associated with {D}edekind sums}, Acta
  Math. Acad. Sci. Hungar. \textbf{26} (1975), no.~3-4, 311--319.

\bibitem{chapmandedekind}
Robin Chapman, \emph{Reciprocity laws for generalized higher dimensional
  {D}edekind sums}, Acta Arith. \textbf{93} (2000), no.~2, 189--199.

\bibitem{dedekind}
Richard Dedekind, \emph{Erl\"auterungen zu den {F}ragmenten xxviii}, Collected
  Works of Bernhard Riemann, Dover Publ., New York, 1953, pp.~466--478.

\bibitem{dieterdedekind}
Ulrich Dieter, \emph{Das {V}erhalten der {K}leinschen {F}unktionen $\log\sigma
  \sb{g,h}(\omega \sb{1}, \omega \sb{2})$ gegen\"uber {M}odultransformationen
  und verallgemeinerte {D}edekindsche {S}ummen}, J. Reine Angew. Math.
  \textbf{201} (1959), 37--70.

\bibitem{ehrhartpolynomial}
Eug{\`e}ne Ehrhart, \emph{Sur les poly\`edres rationnels homoth\'etiques \`a
  {$n$}\ dimensions}, C. R. Acad. Sci. Paris \textbf{254} (1962), 616--618.

\bibitem{girstmair3term}
Kurt Girstmair, \emph{Some remarks on {R}ademacher's three-term relation},
  Arch. Math. (Basel) \textbf{73} (1999), no.~3, 205--207.

\bibitem{hirzebruchzagier}
Friedrich Hirzebruch and Don Zagier, \emph{The {A}tiyah-{S}inger {T}heorem and
  {E}lementary {N}umber {T}heory}, Publish or Perish Inc., Boston, Mass., 1974.

\bibitem{knuth}
Donald~E. Knuth, \emph{The {A}rt of {C}omputer {P}rogramming. {V}ol. 2}, second
  ed., Addison-Wesley Publishing Co., Reading, Mass., 1981.

\bibitem{meyerdedekind}
Curt Meyer, \emph{\"{U}ber einige {A}nwendungen {D}edekindscher {S}ummen}, J.
  Reine Angew. Math. \textbf{198} (1957), 143--203.

\bibitem{meyersczech}
Werner Meyer and Robert Sczech, \emph{\"{U}ber eine topologische und
  zahlentheoretische {A}nwendung von {H}irzebruchs {S}pitzenaufl\"osung}, Math.
  Ann. \textbf{240} (1979), no.~1, 69--96.

\bibitem{mordell}
Louis~J. Mordell, \emph{Lattice points in a tetrahedron and generalized
  {D}edekind sums}, J. Indian Math. Soc. (N.S.) \textbf{15} (1951), 41--46.

\bibitem{pommersheim}
James~E. Pommersheim, \emph{Toric varieties, lattice points and {D}edekind
  sums}, Math. Ann. \textbf{295} (1993), no.~1, 1--24.

\bibitem{rademacherthreeterm}
Hans Rademacher, \emph{Generalization of the reciprocity formula for {D}edekind
  sums}, Duke Math. J. \textbf{21} (1954), 391--397.

\bibitem{rademachersummen}
\bysame, \emph{Zur {T}heorie der {D}edekindschen {S}ummen}, Math. Z.
  \textbf{63} (1956), 445--463.

\bibitem{rademacherdedekind}
\bysame, \emph{Some remarks on certain generalized {D}edekind sums}, Acta
  Arith. \textbf{9} (1964), 97--105.

\bibitem{solomondedsum}
David Solomon, \emph{Algebraic properties of {S}hintani's generating functions:
  {D}edekind sums and cocycles on {${\rm PGL}\sb 2({\bf Q})$}}, Compositio
  Math. \textbf{112} (1998), no.~3, 333--362.

\end{thebibliography}

\def\cprime{$'$} \def\cprime{$'$}
\providecommand{\bysame}{\leavevmode\hbox to3em{\hrulefill}\thinspace}
\providecommand{\MR}{\relax\ifhmode\unskip\space\fi MR }
\providecommand{\MRhref}[2]{%
  \href{http://www.ams.org/mathscinet-getitem?mr=#1}{#2}
}
\providecommand{\href}[2]{#2}

\end{document}